\documentclass[runningheads]{llncs}
\usepackage{cancel}
\usepackage{graphicx}
\usepackage{adjustbox}
\usepackage{overpic}

\usepackage{algorithm}
\usepackage{algorithmic}

\usepackage{mathrsfs}  
\usepackage{sidecap}
\usepackage{bm}
\usepackage{amsfonts}
\usepackage{amssymb}
\usepackage{amsmath}
\usepackage{comment}
\usepackage{xcolor}
\usepackage{relsize}
\usepackage{multirow}
\usepackage{multicol}
\usepackage{diagbox}
\usepackage{tikz}
\usepackage{overpic}
\newcommand{\R}{{\mathbb R}}

\newcommand{\ve}[1]{\mathbf{#1}}
\DeclareMathOperator*{\argmin}{arg\,min}

\def\x{\ve{x}}
\def\y{\ve{y}}

\def\D{\ve{D}}

\newcommand{\circledstar}{%
  \tikz[baseline=(char.base)]{
    \node[shape=circle, draw, inner sep=0.1pt ] (char) {$\star$};
  }%
}

\newtheorem{Proposition}[definition]{Proposition}
\def\mp #1{{\color{black}#1}}

\usepackage{cleveref}
\usepackage{autonum}
\begin{document}
\title{Whiteness-based bilevel estimation of weighted TV parameter maps for image denoising}
\titlerunning{Whiteness-based estimation of weighted TV parameter maps}
%
%
\author{ Monica Pragliola\inst{1} \and Luca Calatroni\inst{2} \and
Alessandro Lanza\inst{3}}

\authorrunning{M. Pragliola, L. Calatroni, A. Lanza}

\institute{Department of Mathematics and Applications, University of Naples Federico II, Naples, Italy  \\\email{monica.pragliola@unina.it} \\
\and
MaLGa Center, DIBRIS, Università di Genova, Istituto Italiano di Tecnologia, Genoa, Italy \\
\email{luca.calatroni@unige.it}\\
\and Department of Mathematics, University of Bologna, Italy\\
\email{alessandro.lanza2@unibo.it}
}
\maketitle              

\begin{abstract}
We consider a bilevel optimisation strategy based on normalised residual whiteness loss for estimating the weighted total variation parameter maps for denoising images corrupted by additive white Gaussian noise. Compared to supervised and semi-supervised approaches relying on prior knowledge of (approximate) reference data and/or information on the noise magnitude, the proposal is fully unsupervised. To avoid noise overfitting
an early stopping strategy is used, relying on simple statistics of optimal performances on a set of natural images.  
Numerical results comparing the supervised/unsupervised procedures for scalar/pixel-dependent \mbox{parameter maps are shown.}

\end{abstract}

\section{Introduction} \label{sec:introduction}

Variational methods are a popular paradigm to overcome the difficulties related to inverse imaging problems: they seek an estimate of the unknown image by minimizing a functional encoding information on the noise and prior beliefs on the image. For image denoising, denoting by $\y\in\R^n$ the (vectorised) data corrupted by Additive White Gaussian Noise (AWGN), variational methods solve
\begin{equation}\label{eq:var_model}
    \x^*(\bm{\lambda})\in\argmin_{\x\in\R^n} \left\{\frac{1}{2}\|\x-\y\|_2^2+R(\x;\bm{\lambda})\right\}\,,
\end{equation}
where $\x^*$ depends on a certain number of parameters $\bm{\lambda}\in\Lambda$ possibly involved in the so-called regularisation term $R$. In general, flexible models are characterised by a large number of parameters for which robust selection strategies need to be introduced. How to properly select parameters in variational models is an issue that has attracted the interest of many researchers working on inverse problems.

Strategies for the selection of $\bm{\lambda}$ maps in variational models can be divided into supervised, semi-supervised and unsupervised. Supervised methods rely on the existence of datasets of groundtruth and corrupted data that allow to customise the regularisation on a specific class of images \cite{Kunisch2013}. Semi-supervised methods 
exploit the availability of partial information, such as the noise statistics - \mbox{see, e.g. the} Discrepancy Principle (DP) \cite{Morozov1966} or its bilevel counterpart \cite{Fehrenbach2015}, and the Stein’s Unbiased Risk
Estimator (SURE), all requiring to know the noise level. Unsupervised methods solely rely on quantities depending on the data $\y$. Empirical approaches for selecting a single regularisation parameter $\bm{\lambda}=\lambda$, such as the L-curve \cite{LC1} or the Generalised Cross Validation \cite{gcv}, belong to this class, which also includes strategies detecting a general number of parameters - see the recent Unknown Noise level SURE (UNSURE)  \cite{tachella2024unsureunknownnoiselevel} and references therein.

Among the unsupervised approaches, in the last years attention has been given to methods observing that a good reconstruction is such that the residual image $\bm{r}^*(\bm{\lambda}):=\y-\x^*(\bm{\lambda})$ is close to the noise realisation degrading the data in terms of its spectrum; this can be done, e.g., enforcing whiteness of the residual, whenever applicable \cite{MF_whiteness}. 
In \cite{etna}, the authors introduced the Residual Whiteness Principle (RWP) in presence of AWGN for setting a scalar regularisation parameter $\lambda$. In \cite{eusipco}, for the same problem, the RWP has been formulated as a bilevel task, where the upper loss measures the whiteness of the residual, and compared to other supervised or semi-supervised losses. Recently, a general bilevel formulation of the RWP, called Generalised Whiteness Principle (GWP), applicable for the estimation of a general number of parameters under the corruption of non white but whitenable noises, has been proposed \cite{COAP}. Nonetheless, the GWP has been only applied to set a small number of parameters and the bilevel task is solved approximately due to the presence of non-smooth lower cost functionals.

In this work, our goal is to extend the bilevel framework introduced in \cite{eusipco} for the estimation of a single parameter $\lambda$ to the case in which the variational model involves a large number of unknown parameters $\bm{\lambda}$. We are exploring the potential of the whiteness loss while mitigating some of the drawbacks that the GWP formulated in \cite{COAP} brings along when applied to such challenging scenario. 
Here, our interest is focused on the selection of the parameters defining a space-variant regularisation term \cite{SIAMrev2021} for the denoising problem under the action of AWGN; more specifically, we are considering the following variational model
\begin{equation}\label{eq:wtvl2}
   \!\!\!\!      \! \x^*(\bm{\lambda}){\in}\argmin_{\x\in\R^n} \left\{\frac{1}{2}\|\x-\y\|_2^2{+}\mathrm{WTV}(\x;\bm{\lambda})\right\}\,,\,\mathrm{WTV}(\x;\bm{\lambda}) {=} \sum_{i=1}^n\lambda_i\|(\bm{\mathrm{D}}\x)_i\|_2\,,
\end{equation}
where the Weighted Total Variation (WTV) term is a weighted version of the TV regulariser \cite{ROF}, and allows to tune the amount of regularisation pixelwise. Numerical tests show that the whiteness principle returns informative $\bm{\lambda}$-maps, and suggest that its employment in highly-parametrised settings could be beneficial.

\section{Bilevel estimation of WTV parameter maps}  \label{sec:bilevel_learn}

We consider a bilevel optimisation approach to compute the optimal parameters $\widehat{\bm{\lambda}}$ for the WTV regulariser in \eqref{eq:wtvl2} given a noisy image $\y\in\R^n$, i.e.~we consider
\begin{align}
    \widehat{\bm{\lambda}}\in & \argmin_{\bm{\lambda}\,\in\,\Lambda}~\mathcal{Q}(\x^*(\bm{\lambda})) :=\frac{1}{2}\| \bm{\rho}(\x^*(\bm{\lambda})) \|_2^2
    \label{eq:bilevel_sup} \\ 
& \text{s.t.}\quad\x^*(\bm{\lambda})\in\argmin_{\x\in\R^n}~\frac{1}{2}\|\x-\y\|_2^2 + \textsc{WTV}(\x;\bm{\lambda}),  \label{eq:bilevel_inf}
\end{align}
where $\mathcal{Q}$ is a metric assessing the quality of the reconstruction $\x^*(\bm{\lambda})$ depending on either the availability of reference data and/or some knowledge on the noise/problem statistics, with $\bm{\rho}(\cdot)$ being chosen accordingly.
Note that although several approximation techniques have been considered to address bilevel problems with non-smooth lower-level constraints (see, e.g., \cite{Ochs2016}), we restrict in the following to the more simple case where in (\ref{eq:bilevel_sup})-(\ref{eq:bilevel_inf}) smoothness is assumed for the lower-level variational constraint, so as to apply implicit differentiation w.r.t.~$\bm{\lambda}$ by means of the implicit function theorem, see, e.g.~\cite{Kunisch2013,Haber_2003}. We will thus consider a smoothed version of \eqref{eq:wtvl2} depending on a smoothing parameter $0<\epsilon \ll 1$,
\begin{equation}  \label{eq:WTV_smoothed}
\textsc{WTV}_\epsilon(\x;\bm{\lambda})  = \sum_{i=1}^n \lambda_i ~h_\epsilon( (\D\x)_i)\,,\quad h_{\epsilon}(\bm{v}) = \begin{cases}
    \frac{3}{4\epsilon}\|\bm{v}
\|_2^2 - \frac{1}{8\epsilon^3}\|\bm{v}\|_2^4\quad&\quad\text{if }\|\bm{v}\|_2<\epsilon\\
\|\bm{v}\|_2 - \frac{3\epsilon}{8}\quad&\quad\text{if }\|\bm{v}\|_2\geq\epsilon\end{cases},
\end{equation}
with $h_\epsilon:\R^2\to \R_{\geq 0}$ a $C^2$-Huber smoothing function. The smoothed functional
\begin{equation} \label{eq:smooth_Feps}
F_\epsilon(\x;\y,\bm{\lambda}) := \frac{1}{2}\|\x-\y\|_2^2 + \textsc{WTV}_\epsilon(\x;\bm{\lambda})
\end{equation}
is thus $C^2$ and $1$-strongly convex. We report in the following its main properties whose proof follows easily from \cite[Proposition III.1]{eusipco}. In the following proposition, $\odot$ denotes element-wise multiplication while $\bm{\Lambda} := [\bm{\lambda}\,;\bm{\lambda}]\in \R^{2n}_{>0}$.

\begin{Proposition}
The functional $F_{\epsilon}$ defined in \eqref{eq:smooth_Feps}
is twice continuously differentiable and $1$-strongly convex on $\R^n$,  hence it admits a unique minimiser. Its gradient $\boldsymbol{\nabla}_{\x} F_{\epsilon} \in \R^n$ and  Hessian $\boldsymbol{\nabla}^2_{\x}  F_{\epsilon} \in \R^{n \times n}$ are given by
\begin{eqnarray}
\boldsymbol{\nabla}_{\x}\, F_{\epsilon}(\x;\y,\bm{\lambda})
&=&
\left( \x - \y \right)
+
\bm{\lambda}\odot\: \bm{\mathrm{D}}^{\mathrm{T}} \, 
\boldsymbol{\nabla}_{\x} H_{\epsilon}(\bm{\mathrm{D}}\x),
\nonumber\\
\boldsymbol{\nabla}^2_{\x}\, F_{\epsilon}(\x;\y,\bm{\lambda})
&=&
\bm{\mathrm{I}}_{n}
+
\D^{\mathrm{T}}\,\mathrm{diag}(\bm{\Lambda})  
\boldsymbol{\nabla}^2  H_{\epsilon}
(\bm{\mathrm{D}}\x) \, \D,
\label{eq:J_Hess}
\end{eqnarray}
where $\boldsymbol{\nabla}  H_{\epsilon}$ and $\boldsymbol{\nabla}^2  H_{\epsilon}$ denote the vector (respectively, matrix) of first-(respectively, second-)order derivatives of $H_\epsilon(\bm{z}):= \sum_{i=1}^n h_\epsilon(z_i)$ w.r.t.~$\bm{z}=\D\x$ and $\mathrm{diag}(\bm{\Lambda})\in\mathbb{R}^{2n\times 2n}_{>0}$ is the diagonal matrix having $\bm{\Lambda}$ on the diagonal.
By defining $\lambda_{\textsc{max}}:= \| \bm{\lambda}\|_{\infty}$, the gradient $\boldsymbol{\nabla}_{\x} F_{\epsilon}$ is thus $L_{\lambda_{\textsc{max}},\epsilon}$-Lipschitz continuous, 
%
and 
\begin{equation}
\small
L_{\lambda_{\textsc{max}},\epsilon} 
\:{=}\:
\max_{\x \in \R^n} \left\| \boldsymbol{\nabla}^2_{\x} F_{\epsilon}(\x;\y,\bm{\lambda}) \right\|_2
\:{\leq}\;
1 + 
\frac{12 \, \lambda_{\textsc{max}}}{\epsilon}
\;{=:}\;\,
\overline{L}_{\lambda_{\textsc{max}},\epsilon} .
\label{eq:L_bound}
\end{equation}
\end{Proposition}
In the context of supervised approaches, the upper loss may be chosen as (half) the mean squared error (MSE) between the reference image $\overline{\x}$ and the reconstruction $\x^*(\bm{\lambda})$ \cite{Kunisch2013,TuomoJMIV2017}, so as to maximise the Signal-to-Noise-Ratio (SNR):
\begin{equation}  \label{eq:supervised_loss}
\mathcal{Q}_{\text{MSE}}(\x^*(\bm{\lambda});\overline{\x}) = \frac{1}{2}\| \bm{\rho}(\x^*(\bm{\lambda})) \|_2^2\,,\quad \bm{\rho}(\x^*(\bm{\lambda})) := \frac{1}{\sqrt{n}}( \x^*(\bm{\lambda})-\overline{\x})\,.
\end{equation}

In our settings, the MSE loss is employed for the computation of benchmark results but its adoption is unfeasible in practice. Assuming that no information is available on the noise level and on the reference image, a natural choice is to consider the bilevel counterpart of the RWP, formally introduced in \cite{eusipco}, that amounts to maximise the whiteness of the residual image $\bm{r}^*(\bm{\lambda})$ via minimisation of (half) the squared norm of its \emph{sample normalised auto-correlation} $\bm{\gamma}(\bm{r}^*(\bm{\lambda}))$:
\begin{equation} \label{eq:whiteness}
\mathcal{Q}_{\text{white}}(\x^*(\bm{\lambda});\y) {=} \frac{1}{2}\| \bm{\rho}(\x^*(\bm{\lambda})) \|_2^2\,,\; \bm{\rho}(\x^*(\bm{\lambda})) {:=}  \boldsymbol{\gamma}(\bm{r}^*(\bm{\lambda})){=}\frac{\bm{r}^*(\bm{\lambda})~ \circledstar~  \bm{r}^*(\bm{\lambda})}{\|\bm{r}^*(\bm{\lambda})\|_2^2}\,, 
\end{equation}
where the normalisation term  eliminates the dependence on $\sigma$, and where
$\x_1~\circledstar~\x_2$ denotes the discrete circular cross-correlation between the matrices $\bm{X}_1, \bm{X}_2\in\R^{n_1\times n_2}$ such that $\text{vec}(\bm{X}_d) = \x_d\in\R^n$, $d=1,2$ and $n_1n_2=n$, which is defined for $(j_1,j_2)\in \left\{0,\ldots,n_1-1\right\}\times\left\{0,\ldots,n_2-1\right\}$  by
\[
(\bm{X}_1~\circledstar~\bm{X}_2)_{j_1,j_2} = \sum_{k_1=0}^{n_1-1} \sum_{k_2=0}^{n_2-1} (\bm{X}_1)_{k_1,k_2}(\bm{X}_2)_{(j_1+k_1)\text{mod}\;n_1,(j_2+k_2)\text{mod}\;n_2}.
\]

%
\noindent The sample normalised auto-correlation $\boldsymbol{\gamma}$ in \eqref{eq:whiteness} is clearly not defined for a null residual $\bm{r}^*(\bm{\lambda}) = \bm{0}_n$ which, however, is obtained only in the limit/degenerate case $\bm{\lambda} = \bm{0}_n$, of no practical interest. Moreover, $\boldsymbol{\gamma}$ in \eqref{eq:whiteness} has the \mbox{properties \cite{etna}}:
\begin{itemize}
    \item[(P$_1$)] $\;\,\boldsymbol{\gamma}_{0,0} = 1$, $\:\boldsymbol{\gamma}_{j_1,j_2} \in [-1,1] \;\: \forall \, (j_1,j_2) \neq (0,0)$;
    \item[(P$_2$)] 
    $\;\,\lim_{n \to \infty} \boldsymbol{\gamma}_{j_1,j_2} = 0 \;\: \forall \, (j_1,j_2) \neq (0,0)$ if $\bm{r}^*(\bm{\lambda})$ is a realisation of AWGN.
\end{itemize}
It follows from (P$_1$) that the metric $\mathcal{Q}_{\text{white}}$ in \eqref{eq:whiteness} has a positive lower bound, $\mathcal{Q}_{\text{white}}(\x^*(\bm{\lambda});\y) \,\;{\geq}\;\, 1/2\,.$
%
%
The asymptotic property (P$_2$), however, does not imply that the lower bound $1/2$ is reached when solving (even exactly) the bilevel problem (\ref{eq:bilevel_sup})-(\ref{eq:bilevel_inf}) with $\mathcal{Q} = \mathcal{Q}_{\text{white}}$ and a generic lower-level optimisation model, for two reasons. Firstly, the space of all possible solutions $\x^*(\bm{\lambda})$ to the lower problem, parameterised by $\bm{\lambda}$, rarely contains the perfect reconstruction $\x^*(\bm{\lambda}) = \overline{\x}$, and it is very unlikely that the space of the associated residual images $\bm{r}^*(\bm{\lambda})$ contains the AWGN realisation actually corrupting $\overline{\x}$. 
Roughly, the probability that the $\bm{r}^*(\bm{\lambda})$-space contains the actual AWGN realisation is higher 
the greater 
the {dimensionality of $\bm{\lambda}$.} Secondly, even when the estimated residual image $\bm{r}^*(\bm{\lambda})$ is a realisation of AWGN - that is 
$\x^*(\bm{\lambda}) = \overline{\x}$ - and even when the size $n = n_1 n_2$ of the image is very large - so that (P$_2$) can be assumed to ``almost hold true'' - 
$\mathcal{Q}_{\text{white}}$ does not even approach the lower bound $1/2$. 
This is the effect, when increasing the size $n$, of the simultaneous decrease of the values $|\boldsymbol{\gamma}_{j_1,j_2}|$ and growth of the number $n$ of $\boldsymbol{\gamma}_{j_1,j_2}$ terms normed in \eqref{eq:whiteness}. 
The results of some Montecarlo simulations (not reported for shortness) suggest that the $\mathcal{Q}_{\text{white}}$ metric of realisations of white Gaussian random noises assume values between $0.9$ \mbox{and $1.1$.}%

\noindent When the $\mathcal{Q}_{\text{white}}$ metric has been coupled with lower optimisation models with few free parameters (one or two), the solution of the bilevel problem was characterised by $\mathcal{Q}_{\text{white}}$ values slightly greater than $1$, and the reconstruction $\x^*(\bm{\lambda})$ was of good quality (see, e.g. Figs. 6.2,6.4 in \cite{etna} and Figs. 4,7 in \cite{COAP}).

However, this is the first time that the $\mathcal{Q}_{\text{white}}$ upper loss is coupled to a highly-parameterised lower model, with the number of free parameters equal to the number of pixels in the image to reconstruct. 
Although on one hand we can expect that the (very large) solution space $\x^*(\bm{\lambda})$ of the lower model contains images very close to the true one $\overline{\x}$ - hence, associated residual images very close to the AWGN realisation - on the other hand there is a real risk that the space contains also bad-quality reconstructions associated to $\mathcal{Q}_{\text{white}}$ values less than the range $0.9-1.1$ characterising white noise realisations. In Section \ref{sec:early}, especially for the latter issue, we will show evidences of the expected behavior.

\section{Optimisation details}  \label{sec:algos}

In this section, we detail the optimisation algorithms used to solve both the lower- and the upper-level problems in (\ref{eq:bilevel_sup})-(\ref{eq:bilevel_inf}). First, we enforce strict positivity of the desired regularisation map by the (element-wise) change of variables $\bm{\lambda} = \exp(\mathbf{\bm{\beta}})$ for $\bm{\beta}\in\R^n$, thus considering the optimisation over the unconstrained variable $\bm{\beta}$ so that $\mathcal{Q}(\x^*(\bm{\lambda})) = \mathcal{Q}(\x^*(\exp(\bm{\beta})))$. Minimising the upper level term in \eqref{eq:bilevel_sup} requires the computation of the gradient of the loss functional w.r.t.~$\bm{\beta}\in\R^n$, which we perform by the chain rule and implicit differentiation as follows:
\begin{align}
&\boldsymbol{\nabla}_{\bm{\beta}}\,\mathcal{Q}(\x^*) =   (\bm{J}_{\bm{\beta}}\, \x^*)^T\notag\,\boldsymbol{\nabla}_{\x} \mathcal{Q}(\x^*) \\
&
   = - \bm{J}_{\bm{\beta}}\,\exp(\bm{\beta})  \bm{J}_{\bm{\lambda}}\,\boldsymbol{\nabla}_{\x}  F_\epsilon(\x^*)  \left( \boldsymbol{\nabla}^2_{\x}~ F_\epsilon(\x^*) \right)^{-1} \boldsymbol{\nabla}_{\x}\mathcal{Q}(\x^*), \label{eq:composite_derivative} 
\end{align}
where we dropped the dependence on $\bm{\beta}$ from $\x^*$ and on $(\y,\exp(\bm{\beta}))$ from $F_\epsilon$. Formula \eqref{eq:composite_derivative} is typically used in the design of gradient- and Newton-type methods for solving \eqref{eq:bilevel_sup};  
the particular expression of $\boldsymbol{\nabla}_{\x}\mathcal{Q}(\cdot)$ depends on the specific choice of the loss function considered, such as $\mathcal{Q}_{\text{MSE}}$ in \eqref{eq:supervised_loss} or $\mathcal{Q}_{\text{White}}$ in \eqref{eq:whiteness}.

\noindent Next, we outline the solvers used to compute at each outer iteration an estimate $\x^*$ of the lower-level problem and to compute a solution $\hat{\bm{\beta}}$ of the bilevel problem.

\subsection{Strongly-convex accelerated gradient descent lower-level solver}  \label{sec:lower_level}

We compute the minimiser $\x^*(\exp(\bm{\beta}))$ of the functional $F_\epsilon$ in \eqref{eq:smooth_Feps} via a strongly-convex variant of Nesterov's accelerated gradient-descent algorithm described, e.g., in \cite{Chambolle-Pock-2016,Calatroni-Chambolle-2019}, see Alg. \ref{algo:Nesterov}, which enjoys linear convergence rate. Due to the $C^2$-regularity of $F_\epsilon$, (quasi) Newton-type techniques could have been used for fast convergence, similarly as in \cite{Kunisch2013,DeLosReyes2015}. This, however, may suffer from high computational costs required for the inversion/approximation of the Hessian along the iterations. On the other hand, the proposed Algorithm is computationally lighter. To get a good approximation of $\x^*$ a small relative tolerance parameter $\texttt{tol}$ is employed and at each iteration $i\geq 1$ of the outer optimisation solver a fixed step-size $1/L_i$,  $L_i=\overline{L}_{\lambda_{\textsc{max}}^i,\epsilon}$, in \eqref{eq:L_bound} depending on the current $\lambda_{\textsc{max}}^i$ is used.

{
\begin{algorithm}[t]\small
\caption{Strongly-convex Nesterov AGD,\,\texttt{SC$\_\text{AGD}$}$(\x_0,L_i,\bm{\lambda}_i,\texttt{tol},\mu=1)$}\label{algo:Nesterov}
\begin{algorithmic}
\STATE \textbf{Initialise:} $\theta_0=1, \tau=1/L_i,~ \x_{-1}=\x_0, t=0$, $\kappa = \mu/L_i$
 \WHILE{$\| \x_{t+1} - \x_{t}\|_2 > \texttt{tol}$}
   \STATE $\theta_{t+1}= \frac{1 - \kappa\theta_t^2+\sqrt{(1-\kappa\theta_t^2)^2+4\theta_t^2}}{2}$; $\quad  \beta_{t+1} = \frac{\theta_t-1}{\theta_{t+1}} \frac{1-\theta_{t+1}\mu\tau}{1-\tau\mu}$
   \STATE $\bm{z}_{t+1} = \x_t + \beta_{t+1}(\x_t-\x_{t-1})$; $\quad\x_{t+1} = \bm{z}_{t+1} - \tau \boldsymbol{\nabla}_{\x} F_{\epsilon}(\bm{z}_{t+1};\y,\bm{\lambda}_i)$
     \STATE $t=t+1$
 \ENDWHILE
 \RETURN $\x^*=\x^*(\bm{\lambda}_i)$
\end{algorithmic}
\end{algorithm}
}
{
\begin{algorithm}[h!]\small
\caption{Gradient-descent for solving \eqref{eq:composite_derivative},\,\texttt{GD$\_$BIL}$(\bm{\beta}^0,\eta,\varepsilon)$}\label{alg:GD_bil}
\begin{algorithmic}
\STATE \textbf{Initialise:} $i=0$, $\x^*(\exp(\bm{\beta}^{-1}))=\y$
 \WHILE{$\| \bm{\beta}^{i+1} - \bm{\beta}^{i}\|_2 > \varepsilon$}
   \STATE compute~$\x^*(\exp(\bm{\beta}^i))$=\texttt{SC}$\_\texttt{AGD}(\x^*(\exp(\bm{\beta}^{i-1})),L_i,\exp(\bm{\beta}^{i}),\texttt{tol},\mu=1)$
        \STATE   compute $\bm{\beta}^{i+1} = \bm{\beta}^{i}-\eta \boldsymbol{\nabla}_{\bm{\beta}}\,\mathcal{Q}(\x^*(\exp(\bm{\beta}^i)))$
     \STATE $i=i+1$
 \ENDWHILE
 \RETURN $\hat{\bm{\lambda}}= \exp(\hat{\bm{\beta}})$
\end{algorithmic}
\end{algorithm}
}

\subsection{Gradient-descent upper-level solver} \label{sec:upper_lev}

We use standard gradient-descent for solving the non-convex bilevel problem \mbox{(\ref{eq:bilevel_sup})-(\ref{eq:bilevel_inf}),} see Alg. \ref{alg:GD_bil}, that takes as input an initial value $\bm{\beta}^0$, a step-size $\eta>0$ and a tolerance $\varepsilon$. While most of the existing approaches make use of quasi-Newton updates for solving efficiently the nested problem \cite{Kunisch2013,TuomoJMIV2017,Fehrenbach2015}, here we focused on a plain first-order scheme. 
Based on \eqref{eq:supervised_loss}, \eqref{eq:whiteness}, we first observe that:
\[
\boldsymbol{\nabla}_{\bm{\beta}} \mathcal{Q}(\x^*) = (\bm{J}_{\bm{\beta}}\,\x^*)^T (\bm{J}_{\bm{x}}\, \bm{\rho})^T \bm{\rho}(\x^*),
\]
where we can proceed as in \eqref{eq:composite_derivative} to compute the Jacobian of $\bm{x}^*$ w.r.t.~$\bm{\beta}$ and of $\bm{\rho}$ w.r.t.~$\bm{x}$ for the specific choice of $\bm{\rho}(\cdot)$ considered. This can be further exploited to design, e.g., Gauss-Newton solvers which could improve speed of convergence. 

\subsection{Early stopping} \label{sec:early}

While in the supervised setting Alg. \ref{alg:GD_bil} can be run until convergence, we observed that running the unsupervised whiteness-based approach until convergence allows in most cases to achieve values of the $\mathcal{Q}_{\text{white}}$ loss {significantly less than the expected interval $0.9-1.1$ and possibly close to the theoretical lower bound $1/2$, which typically corresponds to meaningless parameters $\bm{\lambda}$ and solutions $\x^*(\bm{\lambda})$.

In most of the experiments performed, however, we observed an improvement in the quality metrics during the early iterations of the outer process, which degrade to suboptimal values in later ones. To show such behavior, Alg. \ref{alg:GD_bil} has been run for $3000$ iterations with $\mathcal{Q}=\mathcal{Q}_{\text{White}}$ for recovering image $\#0$ shown  in Fig. \ref{fig:im_tr} starting from a noisy data $\bm{y}$ degraded by AWGN with $\sigma=0.05$ and shown in the left bottom panel of  Fig. \ref{fig:snap}. Along the iterations, we monitored the $\mathcal{Q}$-values as well as the Improved Peak-SNR (IPSNR) and the Improved Structural Similarity Index (ISSIM), measuring the improvement of the output w.r.t. $\y$ in terms of PSNR and SSIM, respectively. The quantities are displayed in the top of Fig. \ref{fig:snap}: note that the whiteness loss takes values below 0.9, while the IPSNR and ISSIM, after an initial increase, start decreasing and stabilise around 0, i.e. no improvement with respect to the data is obtained. In the bottom of Fig. \ref{fig:snap} we show the reconstruction along the iterations of Alg. \ref{alg:GD_bil}: one can notice that at the final iteration the output image resembles the data $\y$. 
\begin{figure}
    \centering{
             {\setlength{\tabcolsep}{0.3pt}
    \begin{tabular}{ccc}
       \includegraphics[width=0.31\linewidth]{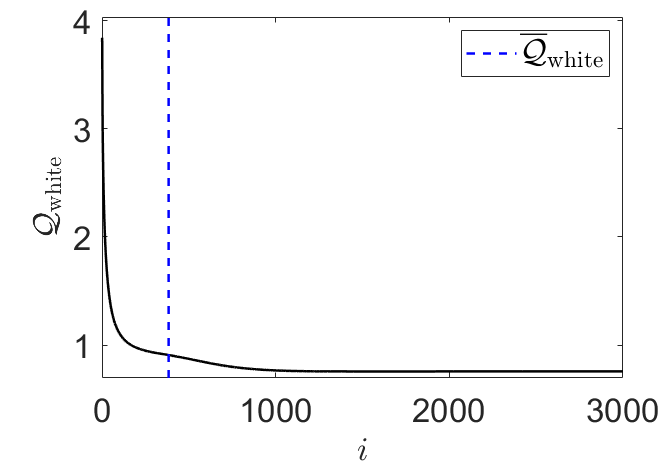}  &   \includegraphics[width=0.3\linewidth]{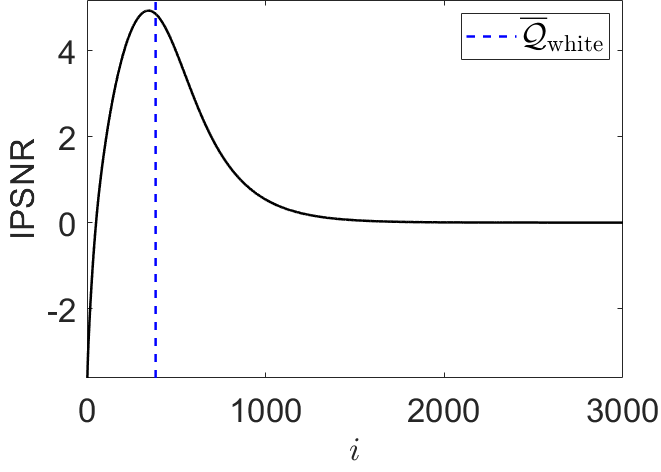} &   \includegraphics[width=0.3\linewidth]{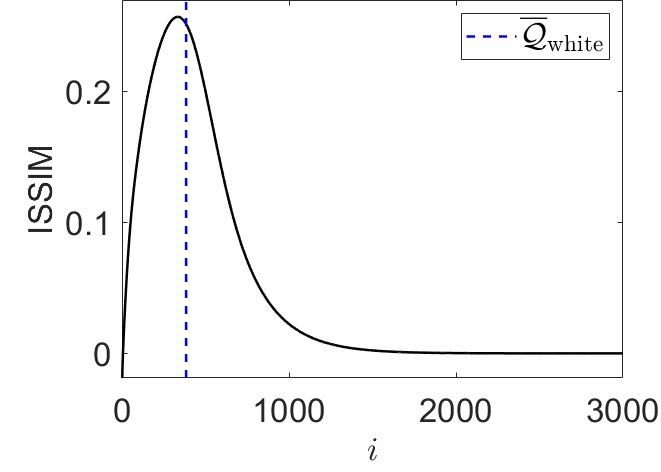}
    \end{tabular}}}
    \centering{
        {\setlength{\tabcolsep}{0.5pt}
    \begin{tabular}{ccccccc}
      ${\bm{y}}$ &$i=1$&$i=100$&$i=300$&$i=600$&$i=1000$&$i=3000$\\
        \includegraphics[width=1.70cm]{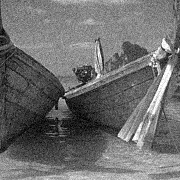}&\includegraphics[width=1.70cm]{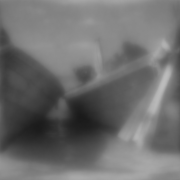} &
        \includegraphics[width=1.70cm]{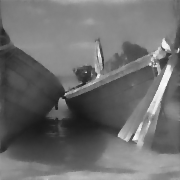} &
        \includegraphics[width=1.70cm]{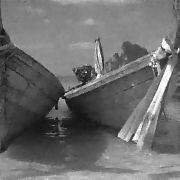} &
        \includegraphics[width=1.70cm]{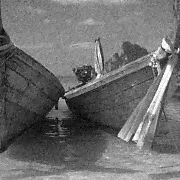}  
        &\includegraphics[width=1.7cm]{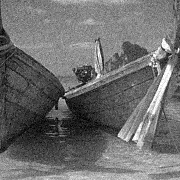}& \includegraphics[width=1.70cm]{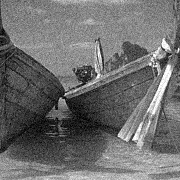} 
        \\
        & \includegraphics[width=1.70cm]{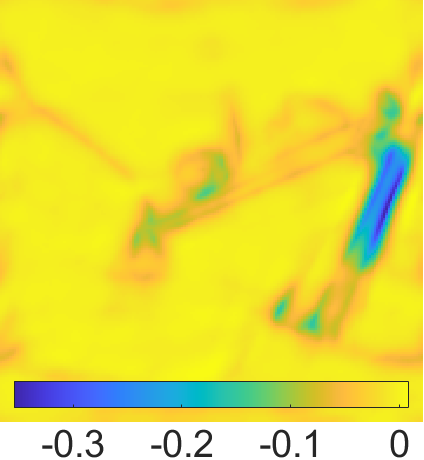} & \includegraphics[width=1.70cm]{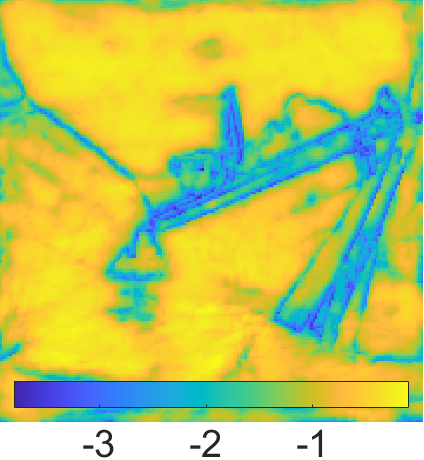} & \includegraphics[width=1.70cm]{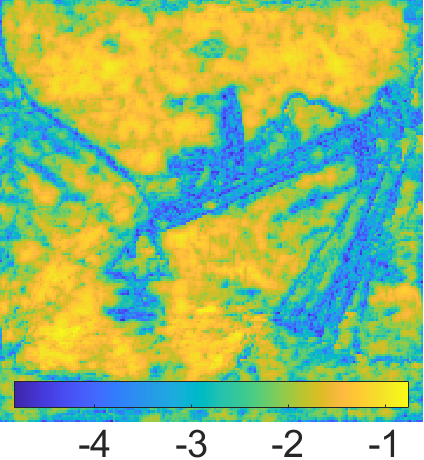} & \includegraphics[width=1.70cm]{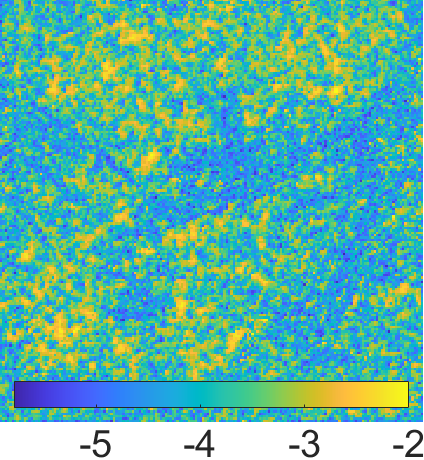} & \includegraphics[width=1.70cm]{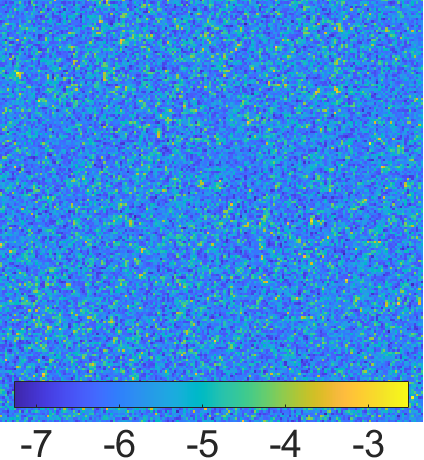} & \includegraphics[width=1.70cm]{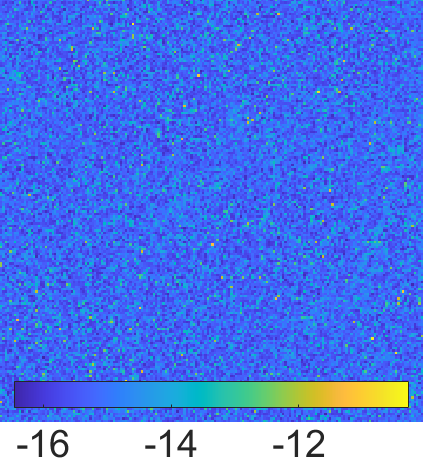}
    \end{tabular}}}
    \caption{Behavior of $\mathcal{Q}$, IPSNR and ISSIM values (top), and of the reconstructions and $\bm{\beta}$-maps (bottom) along the iterations of Alg. \ref{alg:GD_bil} when $\mathcal{Q}=\mathcal{Q}_{\text{white}}$.
    }
    \label{fig:snap}
\end{figure}

By analysing the weights computed along the iterations - see bottom of Fig. \ref{fig:snap} - we observe that at first informative $\bm{\beta}$-values are recovered, which enforce the regularisation on smooth regions of the image while preserving the edges of the structures. However, eventually the $\bm{\beta}$ entries tend to $-\infty$, i.e. the $\bm{\lambda}$-values tend to 0, and the regularisation effect vanishes.


In order to prevent such overfitting, we thus enforced in the optimisation procedure an early stopping regularisation based on statistics of optimal performances estimated off-line. 
{The proposed strategy consists in running the unsupervised procedure until the highest IPSNR value is achieved for a set of {$N=50$} natural images {from the BSD400 repository \cite{BSD400}} corrupted by AWGN with three different standard deviations $\sigma=0.01,0.05,0.1$ and then to compute the values $\mathcal{Q}_{\text{white}}(\x^*_n)$ for all the computed reconstructions $n=1\ldots,N$.} By considering the empirical average of such values, an estimate $\overline{\mathcal{Q}}_{\text{White}}$ of $\mathbb{E}_{\x}(\mathcal{Q}_{\text{white}}(\x))$ is found and can be used to stop Alg. \ref{alg:GD_bil} before overfitting, by replacing the stopping criterion in the \textbf{while} there with a condition checking whether
\begin{equation}  \label{eq:ES}
\| \bm{\beta}^{i+1} - \bm{\beta}^{i}\|_2 > \varepsilon\quad\land\quad \mathcal{Q}_{\text{white}}(\x^*(\exp(\bm{\beta}^{i+1})) \geq \overline{\mathcal{Q}}_{\text{white}}.
\end{equation}
The estimated stopping value is $\overline{\mathcal{Q}}_{\text{White}}=0.9081$. The procedure performs well in practice and does not alter the unsupervised feature of the proposed approach, since inference is performed on an image out-of-distribution. In the previous example, the novel stopping criterion can be employed: the dashed vertical lines in the top panels of Fig. \ref{fig:snap} represent the estimated value $\overline{\mathcal{Q}}_{\text{white}}$ which corresponds to values of quality metrics almost coincident with the highest achievable.

\section{Numerical examples}
\label{sec:test}

In this section, we assess the performance of the whiteness-based procedure
for the bilevel problem \eqref{eq:bilevel_sup}-\eqref{eq:bilevel_inf} when the early stopping criterion outlined in Section \ref{sec:early} is adopted. Our goal is twofold: on one hand we aim to highlight how the unsupervised whiteness loss is capable of estimating a very large number of parameters arising in the variational model while returning high quality results. Moreover, we want to emphasise the potential of the WTV$_{\epsilon}$ regulariser w.r.t the TV$_{\epsilon}$ term, showing how the bilevel problem \eqref{eq:bilevel_sup}-\eqref{eq:bilevel_inf} in presence of the supervised MSE loss can boost the performance of the TV$_{\epsilon}$-based regularisation.
\begin{figure}[b!]
    \centering
    \begin{tabular}{cccc}
       \includegraphics[width=1.90cm]{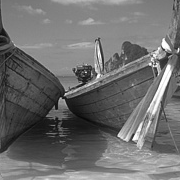}& \includegraphics[width=1.90cm]{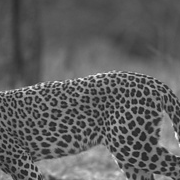}  & \includegraphics[width=1.90cm]{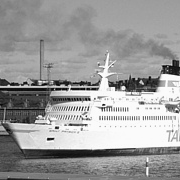}&\includegraphics[width=1.90cm]{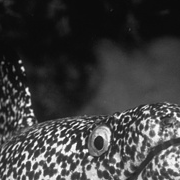}   \\
         
         \#0&\#1&\#2&\#3
    \end{tabular}
    \caption{Original test images of size $180\times 180$ from the BSD400 repository.}
    \label{fig:im_tr}
\end{figure}

\begin{table}[b]\small
    \centering\renewcommand*{\arraystretch}{1.0}
    \setlength{\tabcolsep}{3.5pt}

    \begin{tabular}{c|c|r|r|r|r|r|r|r|r}
    \multicolumn{2}{c}{}&\multicolumn{4}{c}{supervised}&\multicolumn{4}{c}{unsupervised}\\
      \multicolumn{2}{c}{}&\multicolumn{2}{c}{TV}&\multicolumn{2}{c|}{WTV}&\multicolumn{2}{c|}{TV}&\multicolumn{2}{c}{WTV}  \\
      Image&$\sigma$&IPSNR&ISSIM&IPSNR&ISSIM&IPSNR&ISSIM&IPSNR&ISSIM\\
      \hline
      \multirow{3}{*}{\#1}&$0.03$&2.848&0.147&\textbf{6.438}&\textbf{0.194}&2.729&0.129&\textbf{3.176}&\textbf{0.176}\\
&$0.06$&4.257&0.266&\textbf{8.354}&\textbf{0.372}&4.210&0.247&\textbf{4.546}&\textbf{0.322}\\
            &$0.09$&5.345& 0.314&\textbf{9.461}&\textbf{0.471}&5.332&0.303&\textbf{5.462}&\textbf{0.370}\\  
              \hline
\multirow{3}{*}{\#2}&$0.03$&2.896&0.129 &\textbf{6.689}&\textbf{0.179}&2.805&0.117&\textbf{3.265}& \textbf{0.150}\\
&$0.06$&4.467& 0.258&\textbf{8.707}&\textbf{0.362}&4.447&0.248&\textbf{4.817}&\textbf{0.286}\\
&$0.09$&5.624&0.319 &\textbf{10.016}&\textbf{0.465}&5.617&0.313&\textbf{5.891}&\textbf{0.346}\\
              \hline
      \multirow{3}{*}{\#3}&$0.03$&3.309&0.193&\textbf{7.471}&\textbf{0.241}&3.289&0.199&\textbf{3.881} &\textbf{0.224}\\
            &$0.06$&4.604&0.331&\textbf{9.063 }&\textbf{0.434}&4.589& 0.343&\textbf{4.951}&\textbf{0.386}\\   &$0.09$&5.554&0.376&\textbf{10.018}&\textbf{0.525}&5.529&0.395&\textbf{5.794}&\textbf{0.439}
    \end{tabular}
    \caption{IPSNR and ISSIM values achieved by the considered models.}
    \label{tab:1}
\end{table}

\begin{figure}
    \centering
    \begin{tabular}{p{0.4cm}ccccc}
       &&\multicolumn{2}{c}{TV}&\multicolumn{2}{c}{WTV}\\
       &$\bm{y}$&supervised&unsupervised&supervised&unsupervised\\
     \raisebox{0.4cm}{\rotatebox{90}{$\sigma=0.03$}}&
      \begin{overpic}[width=1.90cm]{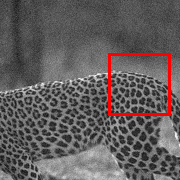} 
\put(-10,-10){\color{red}%
{\includegraphics[scale=0.15]{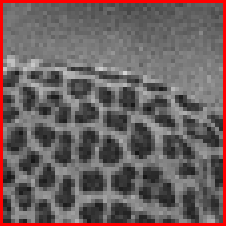} }}
\end{overpic}
      &
       \begin{overpic}[width=1.90cm]{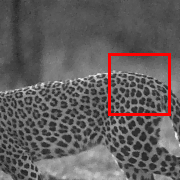} 
\put(-10,-10){\color{red}%
{\includegraphics[scale=0.15]{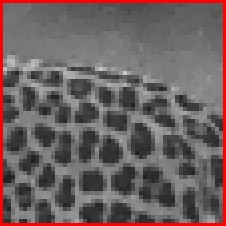} }}
\end{overpic}
    &
\begin{overpic}[width=1.90cm]{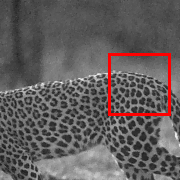} 
\put(-10,-10){\color{red}%
{\includegraphics[scale=0.15]{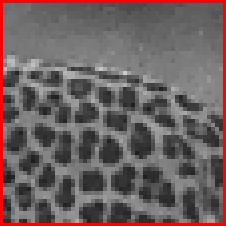} }}
\end{overpic}
    &
    \begin{overpic}[width=1.90cm]{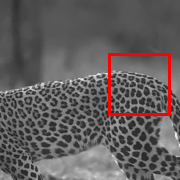} 
\put(-10,-10){\color{red}%
{\includegraphics[scale=0.15]{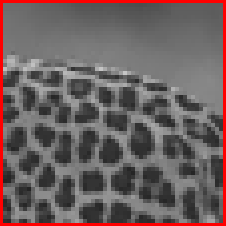} }}
\end{overpic}
&
\begin{overpic}[width=1.90cm]{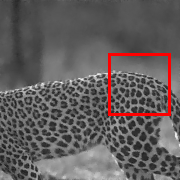} 
\put(-10,-10){\color{red}%
{\includegraphics[scale=0.15]{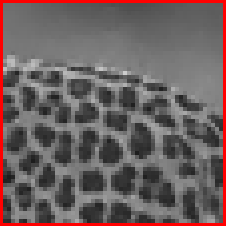} }}
\end{overpic}
   \\
   \\
             \raisebox{0.4cm}{\rotatebox{90}{$\sigma=0.06$}}& \begin{overpic}[width=1.90cm]{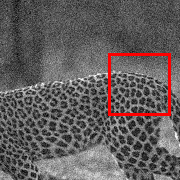} 
\put(-10,-10){\color{red}%
{\includegraphics[scale=0.15]{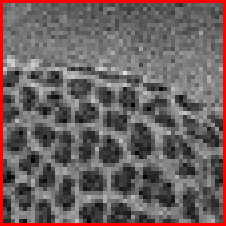} }}
\end{overpic}
      &
       \begin{overpic}[width=1.90cm]{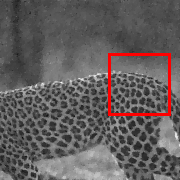} 
\put(-10,-10){\color{red}%
{\includegraphics[scale=0.15]{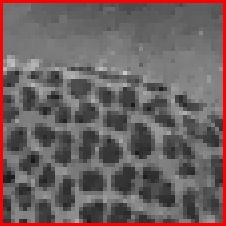} }}
\end{overpic}
    &
\begin{overpic}[width=1.90cm]{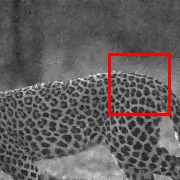} 
\put(-10,-10){\color{red}%
{\includegraphics[scale=0.15]{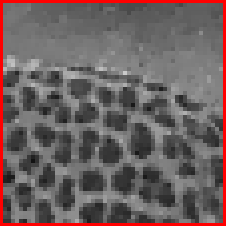} }}
\end{overpic}
    &
    \begin{overpic}[width=1.90cm]{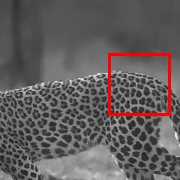} 
\put(-10,-10){\color{red}%
{\includegraphics[scale=0.15]{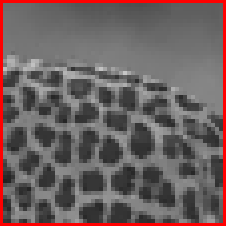} }}
\end{overpic}
&
\begin{overpic}[width=1.90cm]{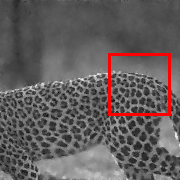} 
\put(-10,-10){\color{red}%
{\includegraphics[scale=0.15]{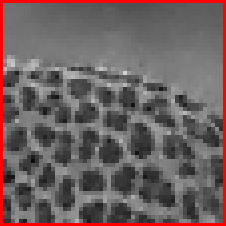} }}
\end{overpic}
        \\
        \\
             \raisebox{0.4cm}{\rotatebox{90}{$\sigma=0.09$}}& \begin{overpic}[width=1.90cm]{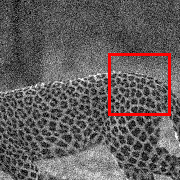} 
\put(-10,-10){\color{red}%
{\includegraphics[scale=0.15]{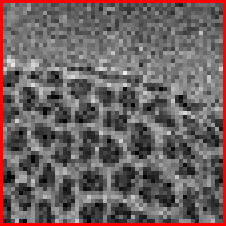} }}
\end{overpic}
      &
       \begin{overpic}[width=1.90cm]{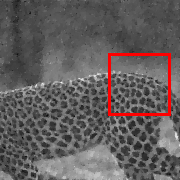} 
\put(-10,-10){\color{red}%
{\includegraphics[scale=0.15]{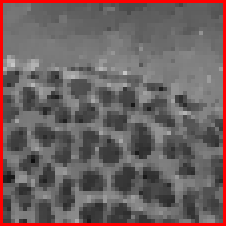} }}
\end{overpic}
    &
\begin{overpic}[width=1.90cm]{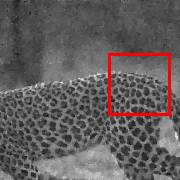} 
\put(-10,-10){\color{red}%
{\includegraphics[scale=0.15]{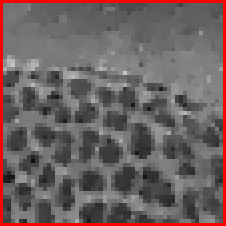} }}
\end{overpic}
    &
    \begin{overpic}[width=1.90cm]{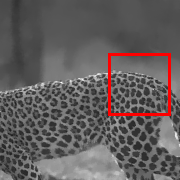} 
\put(-10,-10){\color{red}%
{\includegraphics[scale=0.15]{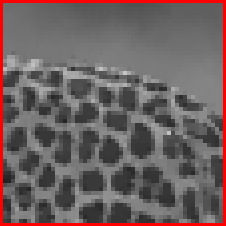} }}
\end{overpic}
&
\begin{overpic}[width=1.90cm]{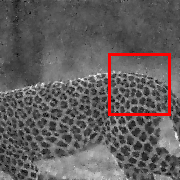} 
\put(-10,-10){\color{red}%
{\includegraphics[scale=0.15]{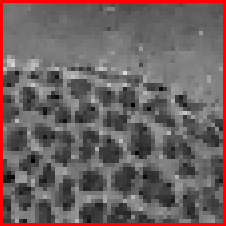} }}
\end{overpic}
    \end{tabular}
    \caption{Output reconstructions for test image \#1.}
    \label{fig:im1}
\end{figure}

Hence, our unsupervised loss is compared to the supervised MSE loss in \eqref{eq:bilevel_sup}. The two losses are also compared in presence of a TV$_{\epsilon}$-regularised lower problem. The four methods are tested on the denoising of test images $\#1$-3 from the BSD400 dataset, shown in Fig. \ref{fig:im_tr}, corrupted by AWGN with different standard deviation $\sigma=0.03,0.06,0.09$ not used in the computation of the stopping value. The smoothing parameters in the smoothed TV$_{\epsilon}$ and WTV$_{\epsilon}$ terms have been selected so as to maximise the performance of the regularised models. More specfically, $\epsilon=10^{-2}$ in the TV$_{\epsilon}$ term and $\epsilon=10^{-1}$ in the WTV$_{\epsilon}$ term. The initial guess $\bm{\beta}_0$ is (a vector of) one(s), while $\eta=1000$ when $\mathcal{Q}=\mathcal{Q}_{\mathrm{white}}$ and $\eta=100$ when $\mathcal{Q}=\mathcal{Q}_{\mathrm{MSE}}$. Also, to avoid degenerate configurations, we impose an upper bound on $\bm{\beta}$, i.e. on $\bm{\lambda}$, $\lambda_i\leq 5$. Finally, $\verb|tol|$ and $\varepsilon$ are set equal to $10^{-6}$.


In Table \ref{tab:1}, for each method we show the IPSNR and ISSIM values obtained, \mp{with in bold the best achieved in the sub-classes of supervised and unsupervised approaches}. The metrics are extremely high when solving \eqref{eq:bilevel_sup}-\eqref{eq:bilevel_inf} with the $\mathcal{Q}_{\text{MSE}}$; also, while the unsupervised method upon the selection of the TV$_{\epsilon}$ term is comparable to the supervised approach, the performance of the WTV$_{\epsilon}$ term in the supervised case does not seem to be reproducible with the unsupervised loss. Nonetheless, in the latter scenario the results outperform the performance of the TV$_{\epsilon}$ term, both in the supervised and unsupervised settings, thus confirming the potentials of the WTV$_{\epsilon}$ term and assessing the robustness of the whiteness loss when moving from a single to a general high number of parameters.

The numeric results reported in Table \ref{tab:1} are confirmed by the output images obtained by the four methods for the different noise levels - see Figs. \ref{fig:im1}-\ref{fig:im3}.

\begin{figure}
    \centering
    \begin{tabular}{p{0.4cm}ccccc}
       &&\multicolumn{2}{c}{TV}&\multicolumn{2}{c}{WTV}\\
       &$\bm{y}$&supervised&unsupervised&supervised&unsupervised\\
    \raisebox{0.4cm}{\rotatebox{90}{$\sigma=0.03$}}&
      \begin{overpic}[width=1.90cm]{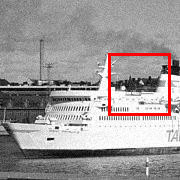} 
\put(-10,-10){\color{red}%
{\includegraphics[scale=0.15]{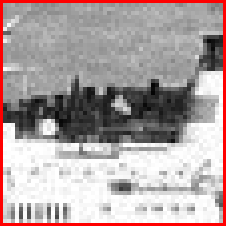} }}
\end{overpic}
      &
       \begin{overpic}[width=1.90cm]{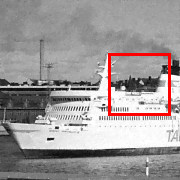} 
\put(-10,-10){\color{red}%
{\includegraphics[scale=0.15]{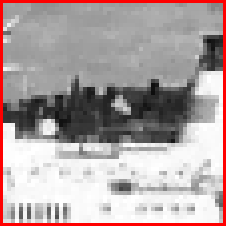} }}
\end{overpic}
    &
\begin{overpic}[width=1.90cm]{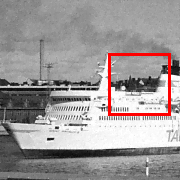} 
\put(-10,-10){\color{red}%
{\includegraphics[scale=0.15]{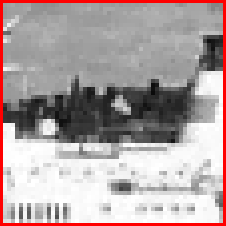} }}
\end{overpic}
    &
    \begin{overpic}[width=1.90cm]{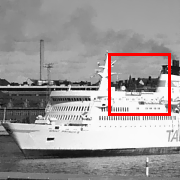} 
\put(-10,-10){\color{red}%
{\includegraphics[scale=0.15]{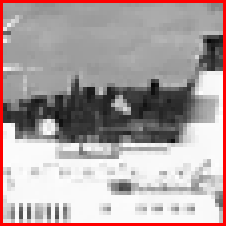} }}
\end{overpic}
&
\begin{overpic}[width=1.90cm]{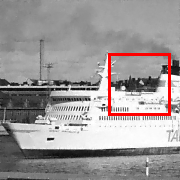} 
\put(-10,-10){\color{red}%
{\includegraphics[scale=0.15]{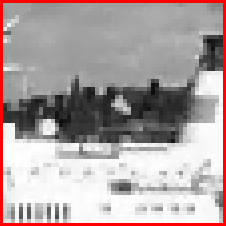} }}
\end{overpic}
   \\
   \\
             \raisebox{0.4cm}{\rotatebox{90}{$\sigma=0.06$}}& \begin{overpic}[width=1.90cm]{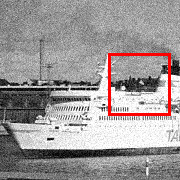} 
\put(-10,-10){\color{red}%
{\includegraphics[scale=0.15]{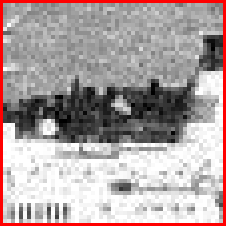} }}
\end{overpic}
      &
       \begin{overpic}[width=1.90cm]{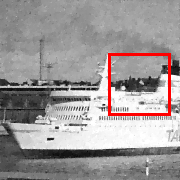} 
\put(-10,-10){\color{red}%
{\includegraphics[scale=0.15]{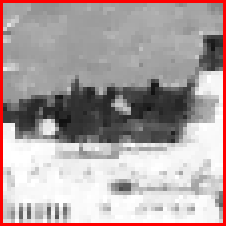} }}
\end{overpic}
    &
\begin{overpic}[width=1.90cm]{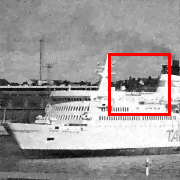} 
\put(-10,-10){\color{red}%
{\includegraphics[scale=0.15]{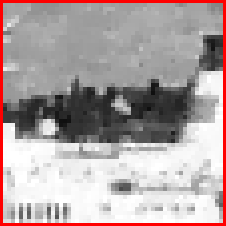} }}
\end{overpic}
    &
    \begin{overpic}[width=1.90cm]{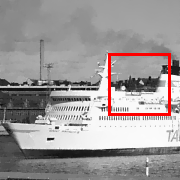} 
\put(-10,-10){\color{red}%
{\includegraphics[scale=0.15]{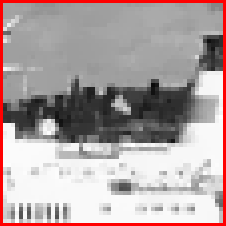} }}
\end{overpic}
&
\begin{overpic}[width=1.90cm]{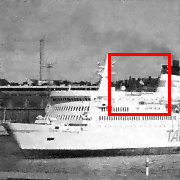} 
\put(-10,-10){\color{red}%
{\includegraphics[scale=0.15]{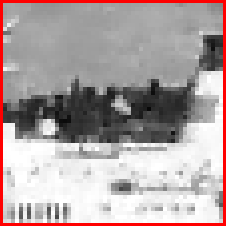} }}
\end{overpic}
        \\
        \\
             \raisebox{0.4cm}{\rotatebox{90}{$\sigma=0.09$}}& \begin{overpic}[width=1.90cm]{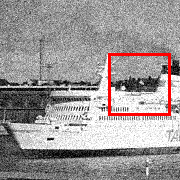} 
\put(-10,-10){\color{red}%
{\includegraphics[scale=0.15]{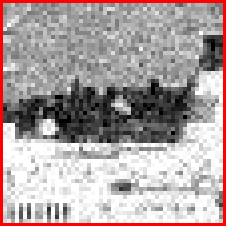} }}
\end{overpic}
      &
       \begin{overpic}[width=1.90cm]{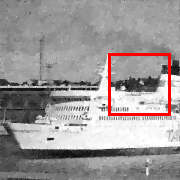} 
\put(-10,-10){\color{red}%
{\includegraphics[scale=0.15]{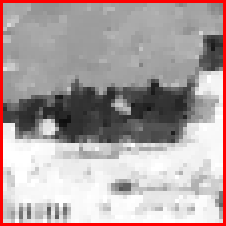} }}
\end{overpic}
    &
\begin{overpic}[width=1.90cm]{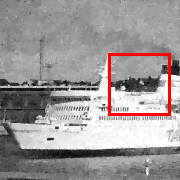} 
\put(-10,-10){\color{red}%
{\includegraphics[scale=0.15]{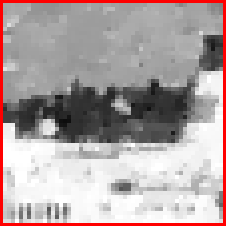} }}
\end{overpic}
    &
    \begin{overpic}[width=1.90cm]{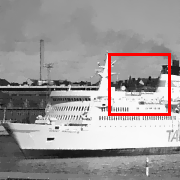} 
\put(-10,-10){\color{red}%
{\includegraphics[scale=0.15]{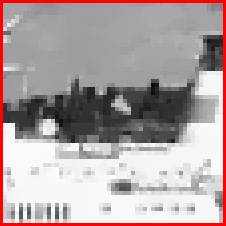} }}
\end{overpic}
&
\begin{overpic}[width=1.90cm]{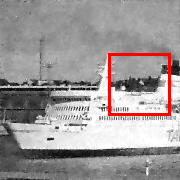} 
\put(-10,-10){\color{red}%
{\includegraphics[scale=0.15]{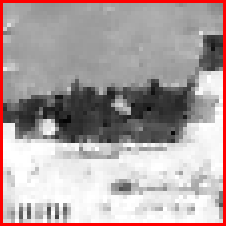} }}
\end{overpic}
    \end{tabular}
    \caption{Output reconstructions for test image \#2.}
  \label{fig:im2}
\end{figure}
Finally, in Fig. \ref{fig:maps} we show the maps of the output $\bm{\beta}$ for the WTV$_{\epsilon}$ regulariser in the supervised and unsupervised scenario; for rendering purposes, a lower bound of $-10$ has been considered in the visualisation of the maps. The supervised and unsupervised approach are equally good at differentiating between smooth regions, characterised by larger values of the parameters, and textured parts of the image, for which smaller values are more suitable. However the scale of the maps in the two cases are significantly different, this possibly being responsible for the loss of quality observed in the unsupervised case.

\begin{figure}
    \centering
    \begin{tabular}{p{0.4cm}ccccc}
       &&\multicolumn{2}{c}{TV}&\multicolumn{2}{c}{WTV}\\
       &$\bm{y}$&supervised&unsupervised&supervised&unsupervised\\
    \raisebox{0.4cm}{\rotatebox{90}{$\sigma=0.03$}}&
      \begin{overpic}[width=1.90cm]{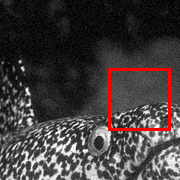} 
\put(-10,-10){\color{red}%
{\includegraphics[scale=0.15]{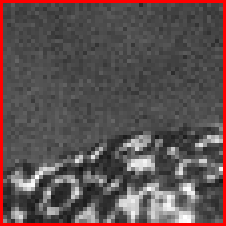} }}
\end{overpic}
      &
       \begin{overpic}[width=1.90cm]{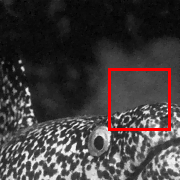} 
\put(-10,-10){\color{red}%
{\includegraphics[scale=0.15]{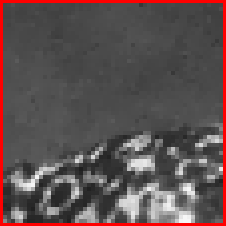} }}
\end{overpic}
    &
\begin{overpic}[width=1.90cm]{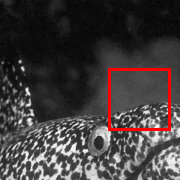} 
\put(-10,-10){\color{red}%
{\includegraphics[scale=0.15]{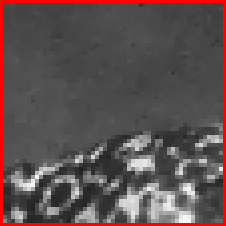} }}
\end{overpic}
    &
    \begin{overpic}[width=1.90cm]{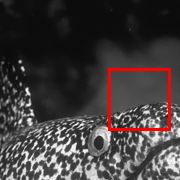} 
\put(-10,-10){\color{red}%
{\includegraphics[scale=0.15]{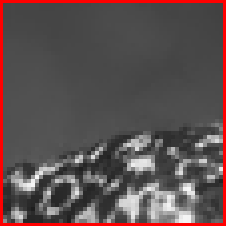} }}
\end{overpic}
&
\begin{overpic}[width=1.90cm]{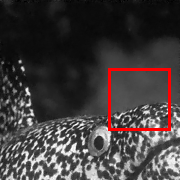} 
\put(-10,-10){\color{red}%
{\includegraphics[scale=0.15]{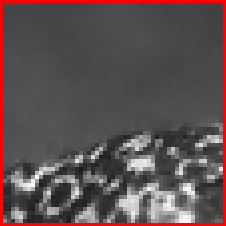} }}
\end{overpic}
   \\
   \\
              \raisebox{0.4cm}{\rotatebox{90}{$\sigma=0.06$}}& \begin{overpic}[width=1.90cm]{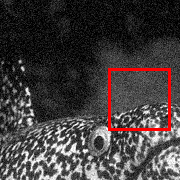} 
\put(-10,-10){\color{red}%
{\includegraphics[scale=0.15]{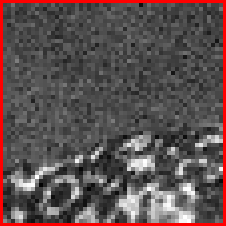} }}
\end{overpic}
      &
       \begin{overpic}[width=1.90cm]{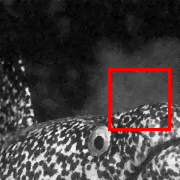} 
\put(-10,-10){\color{red}%
{\includegraphics[scale=0.15]{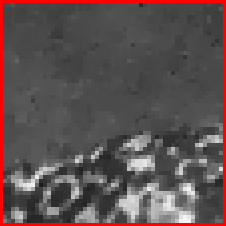} }}
\end{overpic}
    &
\begin{overpic}[width=1.90cm]{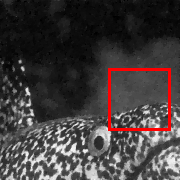} 
\put(-10,-10){\color{red}%
{\includegraphics[scale=0.15]{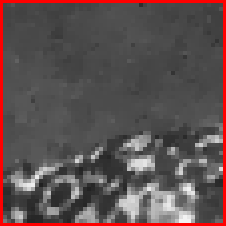} }}
\end{overpic}
    &
    \begin{overpic}[width=1.90cm]{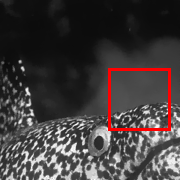} 
\put(-10,-10){\color{red}%
{\includegraphics[scale=0.15]{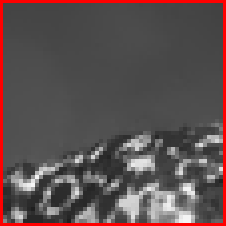} }}
\end{overpic}
&
\begin{overpic}[width=1.90cm]{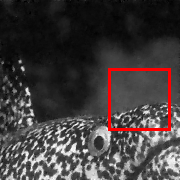} 
\put(-10,-10){\color{red}%
{\includegraphics[scale=0.15]{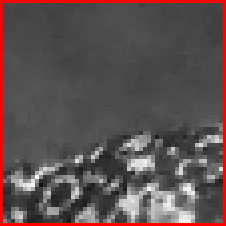} }}
\end{overpic}
        \\
        \\
            \raisebox{0.4cm}{\rotatebox{90}{$\sigma=0.09$}}& \begin{overpic}[width=1.90cm]{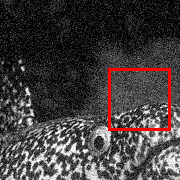} 
\put(-10,-10){\color{red}%
{\includegraphics[scale=0.15]{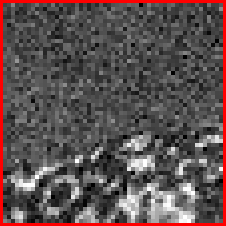} }}
\end{overpic}
      &
       \begin{overpic}[width=1.90cm]{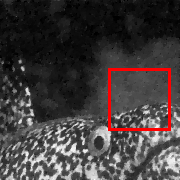} 
\put(-10,-10){\color{red}%
{\includegraphics[scale=0.15]{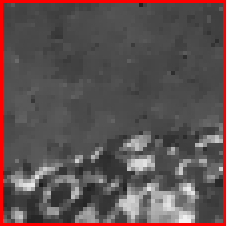} }}
\end{overpic}
    &
\begin{overpic}[width=1.90cm]{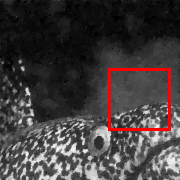} 
\put(-10,-10){\color{red}%
{\includegraphics[scale=0.15]{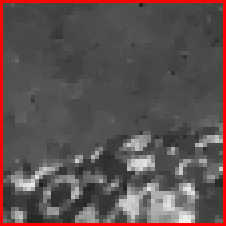} }}
\end{overpic}
    &
    \begin{overpic}[width=1.90cm]{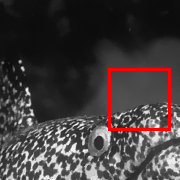} 
\put(-10,-10){\color{red}%
{\includegraphics[scale=0.15]{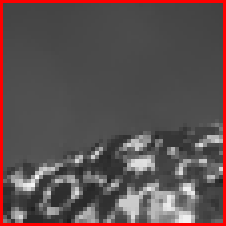} }}
\end{overpic}
&
\begin{overpic}[width=1.90cm]{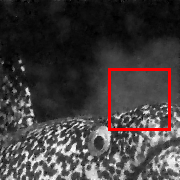} 
\put(-10,-10){\color{red}%
{\includegraphics[scale=0.15]{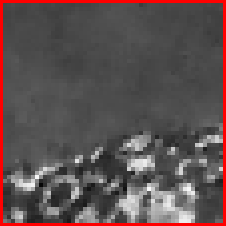} }}
\end{overpic}
    \end{tabular}
  \caption{Output reconstructions for test image \#3.}
  \label{fig:im3}
\end{figure}

\begin{figure}
    \centering
    \setlength{\tabcolsep}{0.5pt}
    \begin{tabular}{ccccccc}
  &\multicolumn{2}{c}{$\sigma=0.03$}&\multicolumn{2}{c}{$\sigma=0.06$}&\multicolumn{2}{c}{$\sigma=0.09$}\\
&supervised&unsupervised&supervised&unsupervised&supervised&unsupervised   \\
\raisebox{1.15cm}{\#1}&\includegraphics[width=0.162\linewidth]{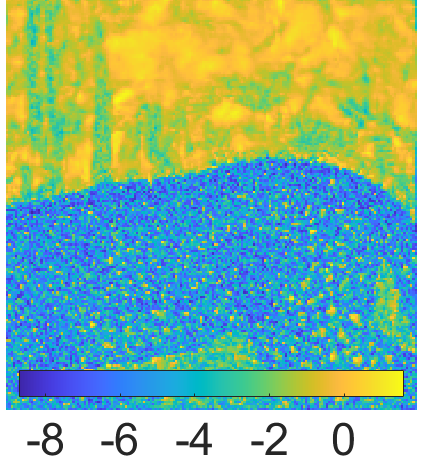}&\includegraphics[width=0.162\linewidth]{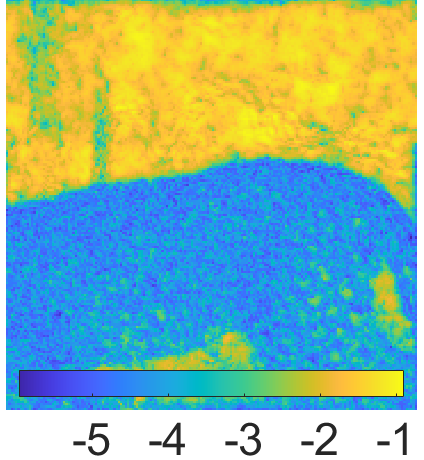}&\includegraphics[width=0.162\linewidth]{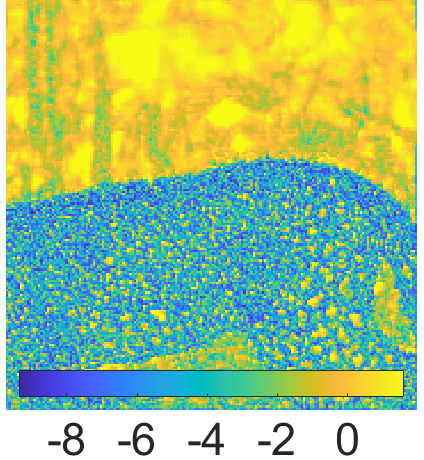}&\includegraphics[width=0.162\linewidth]{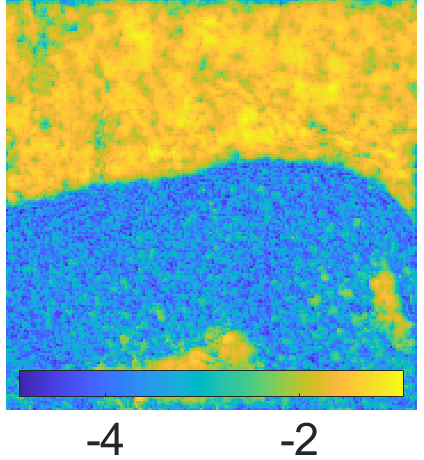}&\includegraphics[width=0.162\linewidth]{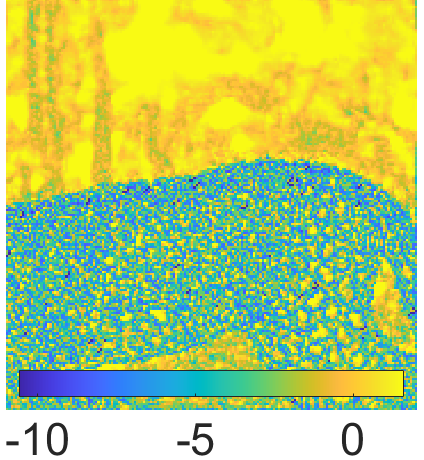}&\includegraphics[width=0.162\linewidth]{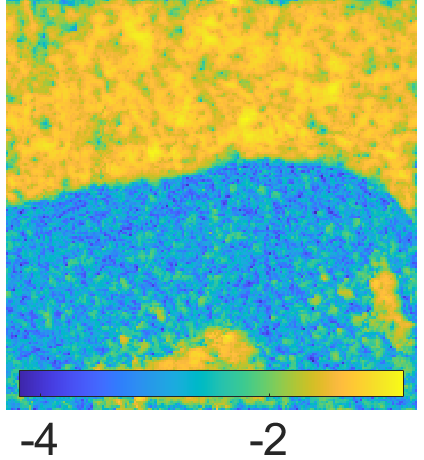}\\

\raisebox{1.15cm}{\#2}&\includegraphics[width=0.162\linewidth]{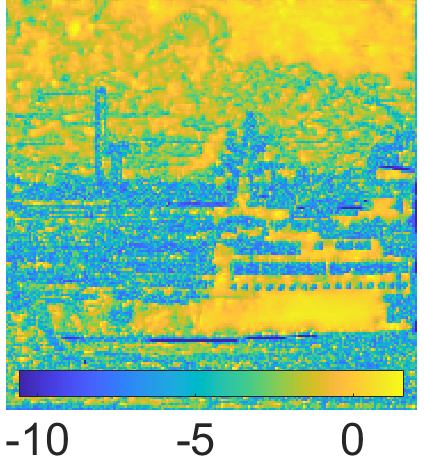}&\includegraphics[width=0.162\linewidth]{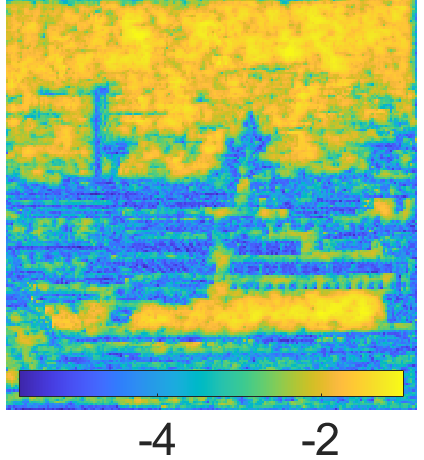}&\includegraphics[width=0.162\linewidth]{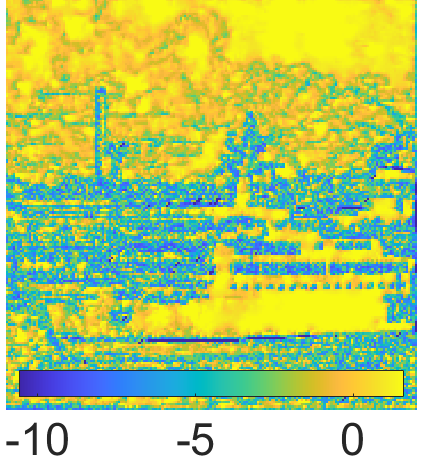}&\includegraphics[width=0.162\linewidth]{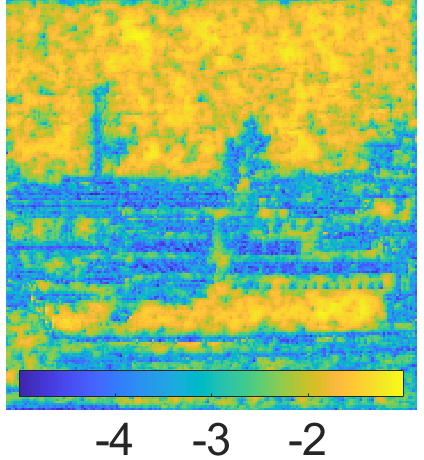}&\includegraphics[width=0.162\linewidth]{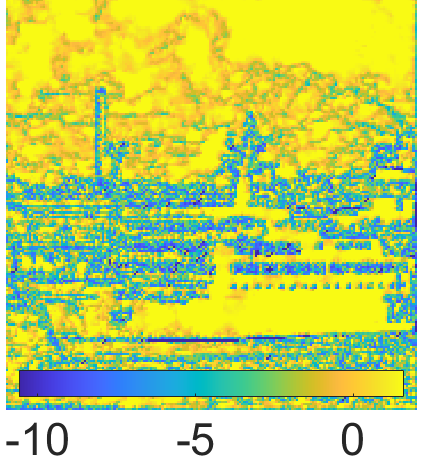}&\includegraphics[width=0.162\linewidth]{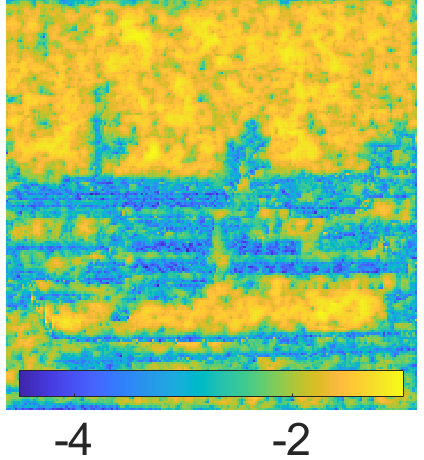}\\

\raisebox{1.15cm}{\#3}&\includegraphics[width=0.162\linewidth]{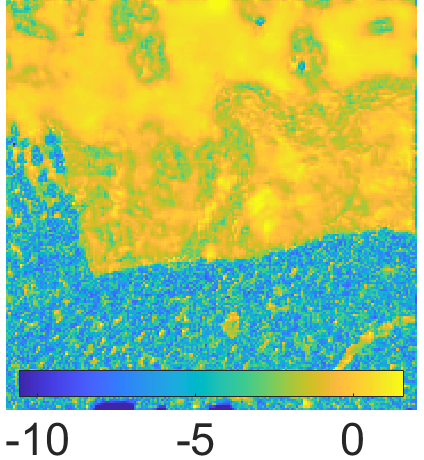}&\includegraphics[width=0.162\linewidth]{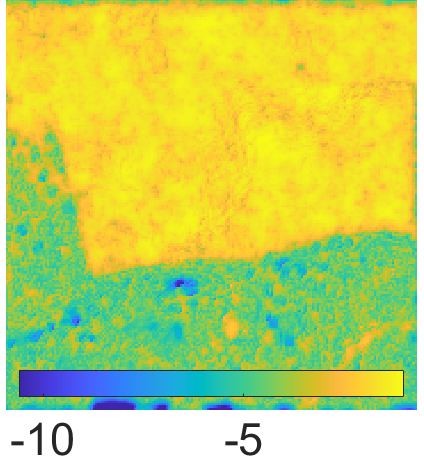}&\includegraphics[width=0.162\linewidth]{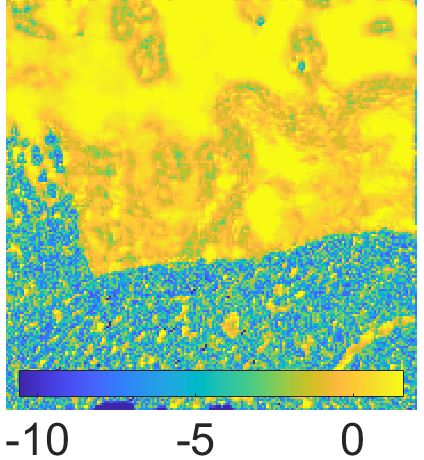}&\includegraphics[width=0.162\linewidth]{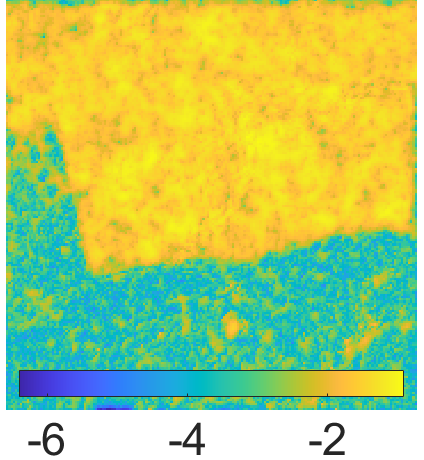}&\includegraphics[width=0.162\linewidth]{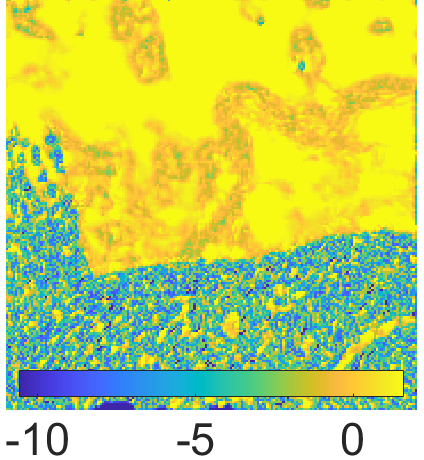}&\includegraphics[width=0.162\linewidth]{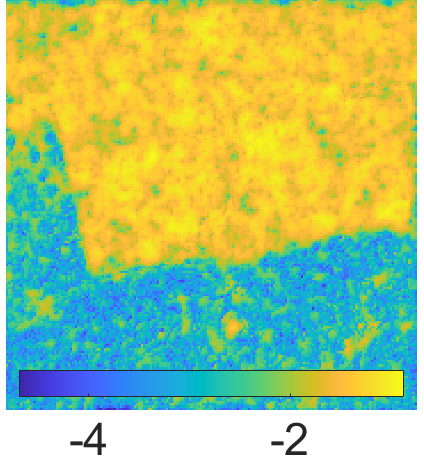}
    \end{tabular}
    \caption{Output parameters for the supervised $\mathcal{Q}_{\mathrm{MSE}}$ and the unsupervised $\mathcal{Q}_{\mathrm{white}}$.}
    \label{fig:maps}
\end{figure}

\mp{The {unsupervised} whiteness-based strategy with the WTV regulariser relies on a training on the BSD400 natural images dataset. To prove the flexibility of the approach, we also perform a training on the COCO dataset of common objects \cite{cocodataset}. From $N=50$ images of the set, we extract a $180\times 180$ patch and consider the same degradation of Section \ref{sec:early}. The learned stopping value is $\widetilde{\mathcal{Q}}_{\text{white}}=0.9190$, very close to the value $\bar{\mathcal{Q}}_{\text{white}}$ employed in the experimentation so far. In Fig. \ref{fig:qs}, for image $\#2$ and the three noise levels, we plot the IPSNRs achieved along a fixed number of iterations of Alg. \ref{alg:GD_bil} for solving the whiteness-based bilevel probelm with the WTV regulariser. The vertical blue and red lines indicate the stopping values returned by the training on BSD400 and COCO, respectively. One can notice that they correspond to very close IPSNRs.}  

\begin{figure}
    \centering
    \begin{tabular}{ccc}
      \includegraphics[width=0.31\textwidth]{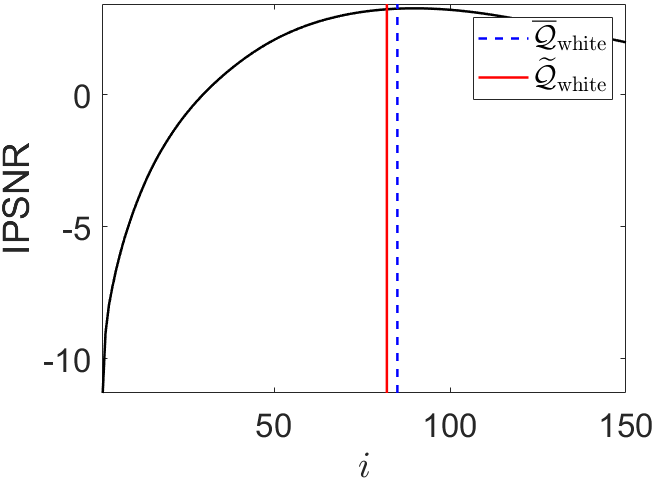}   &  \includegraphics[width=0.31\textwidth]{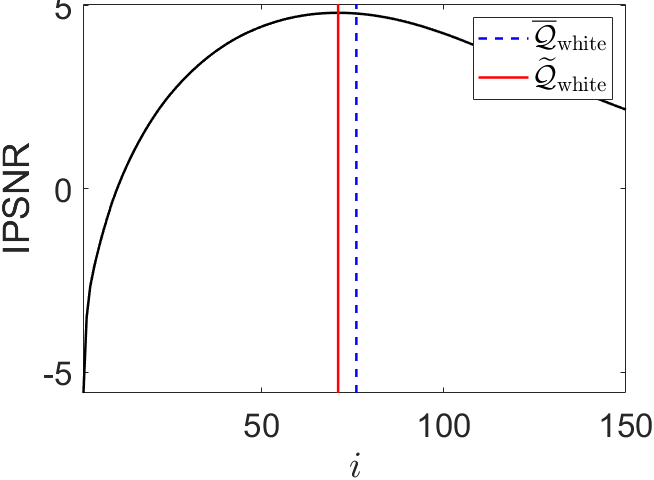}& \includegraphics[width=0.31\textwidth]{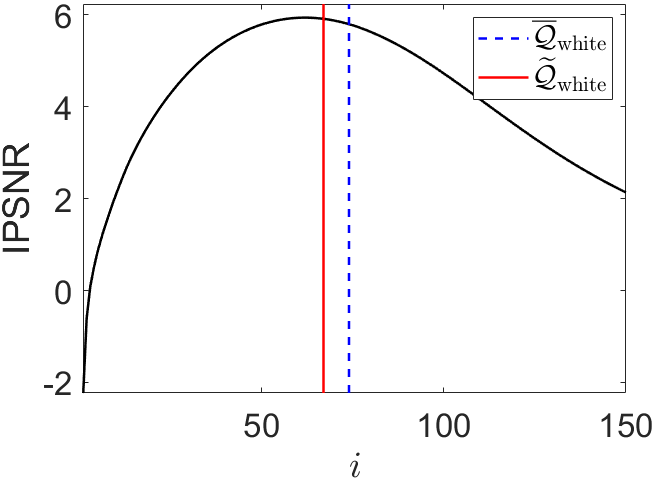}\\
      $\sigma=0.03$&$\sigma=0.06$&$\sigma=0.09$
    \end{tabular}
    \caption{IPSNR along the iterations of Alg. \ref{alg:GD_bil} for different noise levels.}
    \label{fig:qs}
\end{figure}

\mp{Finally, we consider a signal $\x$ corrupted by additive white uniform noise distribution, and solve the bilevel problem in \eqref{eq:bilevel_sup}-\eqref{eq:bilevel_inf} under the same settings detailed above for Gaussian data. The resulting IPSNRs/ISSIMs for image $\#3$ are summarised in Table \ref{tab:2}, and show the same behavior of those reported in Table \ref{tab:1}, suggesting the robustness of the unsupervised whiteness loss also when coupled to lower cost functionals that do not properly model the noise at hand.}

\begin{table}[t]\small
    \centering\renewcommand*{\arraystretch}{1.0}
    \setlength{\tabcolsep}{3.5pt}

    \begin{tabular}{c|c|r|r|r|r|r|r|r|r}
    \multicolumn{2}{c}{}&\multicolumn{4}{c}{supervised}&\multicolumn{4}{c}{unsupervised}\\
      \multicolumn{2}{c}{}&\multicolumn{2}{c}{TV}&\multicolumn{2}{c|}{WTV}&\multicolumn{2}{c|}{TV}&\multicolumn{2}{c}{WTV}  \\
      Image&$\sigma$&IPSNR&ISSIM&IPSNR&ISSIM&IPSNR&ISSIM&IPSNR&ISSIM\\
        \hline
      \multirow{3}{*}{\#3}&$0.03$&3.419&0.199&\textbf{7.377}&\textbf{0.243}&3.354&0.209&\textbf{3.815}&\textbf{0.225}\\
        &$0.06$&4.712&0.342&\textbf{8.900}&\textbf{0.436}&4.629&0.365&\textbf{5.007}&\textbf{0.392}\\
        &$0.09$&5.697&0.387&\textbf{9.887}&\textbf{0.524}&5.585&0.423&\textbf{5.831} &\textbf{0.427}    \end{tabular}
    \caption{IPSNR and ISSIM values achieved by the considered models in case of additive white uniform noise corruption.}
    \label{tab:2}
\end{table}

\section{Conclusions}
\label{sec:concl}

We addressed the estimation of spatially-adaptive weighted TV$_\epsilon$ regularisation maps through a bilevel optimisation framework with an unsupervised upper loss defined in terms of a residual whiteness measure. This approach enables the simultaneous recovery of both the denoised image and the optimal regularisation maps for data corrupted by AWGN.
The approach requires the employment of an early stopping criterion based on simple statistics of optimal performances computed on a set
of natural images, and returns maps morphologically similar to the ones obtained by employing supervised MSE loss, while outperforming the standard (global) TV$_\epsilon$ term, both in supervised and unsupervised settings. Our tests are encouraging in terms of the expected performance of the whiteness loss in highly-parametrised settings. Further work should address the use of the whiteness loss in proper learning settings where, e.g., the $\bm{\lambda}$ map is parametrised as the output of a Convolutional Neural Network (CNN), as in \cite{Kofler2023} where a similar strategy is performed for MRI image reconstruction in a supervised setting.

\section*{Acknowledgments}
This work is based upon work supported by the ANR JCJC project TASKABILE ANR-22-CE48-0010 and by the ERC grant MALIN under the European Union’s Horizon Europe programme (grant agreement No. 101117133). Views and opinions expressed are however those of the author(s) only and do not necessarily reflect those of the European Union or the European Research Council Executive Agency. Neither the European Union nor the granting authority can be held responsible for them. MP acknowledges the GNCS-INdAM funding and the 2022 research program FRA of the University of Naples Federico II.

\bibliographystyle{plain}

\bibliography{refs}

\end{document}